\documentclass{ifacconf}

\usepackage{booktabs}
\usepackage{graphicx}      
\usepackage{natbib}        
\usepackage{algpseudocode}
\usepackage{xfrac}
\usepackage{epic}
\usepackage{amsfonts}
\usepackage{amsmath,amssymb}
\usepackage{textcomp}
\usepackage{bbold}

\makeatletter
\def\endfigure{\end@float}
\def\endtable{\end@float}
\makeatother

\let\ifacconfcaptionwidth\captionwidth
\usepackage[caption=false]{subfig}
\let\captionwidth\ifacconfcaptionwidth

\renewcommand{\min}{\text{min}}

\usepackage{color}
\definecolor{wheat}{rgb}{0.96,0.87,0.70}

\usepackage{calc}

\graphicspath{{figures/}}

\begin{document}

\makeatletter
\renewcommand*\env@matrix[1][*\c@MaxMatrixCols c]{%
  \hskip -\arraycolsep
  \let\@ifnextchar\new@ifnextchar
  \array{#1}}
\makeatother

\begin{frontmatter}

\title{HPIPM: a high-performance quadratic programming framework for model predictive control\thanksref{footnoteinfo}}
\thanks[footnoteinfo]{This research was supported by the German Federal Ministry for Economic Affairs and Energy (BMWi) via eco4wind (0324125B) and DyConPV (0324166B), and by DFG via Research Unit FOR 2401.}

\author[IMTEK]{Gianluca Frison} 
\author[IMTEK]{Moritz Diehl}

\address[IMTEK]{Department of Microsystems Engineering, University of Freiburg, email: \{gianluca.frison, moritz.diehl\} at imtek.uni-freiburg.de}

\begin{abstract}                

This paper introduces HPIPM, a high-performance framework for quadratic programming (QP), designed to provide building blocks to efficiently and reliably solve model predictive control problems.
HPIPM currently supports three QP types, and provides interior point method (IPM) solvers as well (partial) condensing routines.
In particular, the IPM for optimal control QPs is intended to supersede the HPMPC solver, and it largely improves robustness while keeping the focus on speed.
Numerical experiments show that HPIPM reliably solves challenging QPs, and that it outperforms other state-of-the-art solvers in speed.

\end{abstract}

\begin{keyword}
quadratic programming, model predictive control, embedded optimization, software 
\end{keyword}

\end{frontmatter}

\section{Introduction}

The aim of this paper is to introduce the quadratic programming (QP) framework in the high-performance software HPIPM in the context of model predictive control (MPC)~\citep{Maciejowski2002}.
Among other features, the HPIPM software provides an interior point method (IPM) for the solution of linear-quadratic optimal control problems (OCP), which has been used in numerous applications.
In this regard it is intended to supersede the software HPMPC~\citep{Frison2014}, which comprises both hardware-tailored linear algebra (LA) routines and solvers for QPs in the MPC form. 
The LA framework in HPMPC has been completely reimplemented and extended in the BLASFEO software~\citep{Frison2018}, which targets routines operating on both the panel-major matrix format proposed in HPMPC, and the column-major matrix format used in the standard BLAS application programming interface (API). 
Analogously, the QP framework in HPMPC has been completely reimplemented and extended in the HPIPM software, which is the subject of the current paper.

More in detail, HPIPM is an open-source C-coded high-performance framework for QP, and it is designed in a modular fashion with two distinct goals: to be easily extendable to handle new QP types, and to serve as the basis for more generic algorithms such as solvers for nonlinear programming (NLP).
At the time of writing, HPIPM supports three QP types: a dense QP, an OCP QP, and a tree-structured OCP QP.
For each of them, it defines C structures for the QP itself, its solution and its residuals.
These structures serve as arguments to a set of routines factorizing and solving the Karush-Kuhn-Tucker (KKT) system and implementing a family of IPM solvers.
The different IPM solvers are meant to provide different trade-offs between computational speed and robustness, with the fastest but least robust algorithm being analogous to the one implemented in HPMPC.
At their core, the IPM solvers in HPIPM comprise a common set of routines that operate only on vectors and are agnostic of the QP type or even the matrix format employed.
Conversely, the KKT modules, which account for the computationally most intensive operations, are tailored to the specific QP type, and are currently implemented using BLASFEO. 

Alongside these QP types and respective IPM solvers, HPIPM implements routines to convert between them.
Condensing and partial condensing~\citep{Axehill2015} routines convert an OCP QP into a dense QP and an OCP QP with shorter horizon, respectively.
These (partially) condensed QPs can be solved with tailored IPM solvers, and their solution expanded back into the original OCP QP format by means of dedicated routines.
The (partial) condensing and expansion routines can also be employed as stand-alone routines, and coupled with other QP solvers. 
As such, HPIPM provides a complete QP framework.

HPIPM uses encapsulation mechanisms for its C API:
the C structures are not assumed to be accessed directly, but through setters and getters.
This fact has several important consequences: for example, it allows the problem data to be stored in the C structures in a format efficient for the solver 
while setters and getters are designed to be convenient for the user. 
Additionally, this enhances maintainability, as the internal layout of the structures can freely change as long as the setters and getters expose the same API.
Another key feature of HPIPM is its memory management, which eases its embedded use: no internal memory allocation takes place, as HPIPM operates on externally provided chunks of memory. 
On top of the C API, HPIPM provides interfaces to the high-level programming languages Matlab, Octave, Simulink and Python.

At the time of writing, HPIPM is employed as a QP solver in several NLP solvers such as acados~\citep{Verschueren2019}, Control Toolbox~\citep{Giftthaler2018}, MPCC~\citep{Liniger2015} and MATMPC~\citep{Chen2018}.
In particular, in the case of acados it provides the entire QP framework and it is deeply embedded in the source code, in the sense that the members of the OCP QP type are directly accessed when building the QP approximation at each SQP iteration.
Furthermore, the (partial) condensing routines provided by HPIPM are used as pre-processing step in the solution of these QPs by means of other QP solvers such as qpOASES~\citep{Ferreau2014}, qpDUNES~\citep{Frasch2014}, OSQP~\citep{Stellato2017}, OOQP~\citep{Gertz2003}.
In the case of the Control Toolbox, HPIPM is used through its setter-based C API as a QP solver for both constrained and unconstrained linear-quadratic OCPs and it provides the feedback gains needed in the iterative linear quadratic regulator (iLQR) algorithm~\citep{Li2004a}.
In the case of MATMPC, HPIPM is used through its C API within mex sources in Matlab, for both the IPM solver and the (partial) condensing routines.
In the case of MPCC, HPIPM is currently the solver of choice for real-time applications (replacing the originally employed FORCES solver~\citep{Domahidi2012}), 
and it is used through its Matlab interface. 


\section{QP formulations}

HPIPM currently defines three QP types: a dense QP, an OCP QP and tree-structured OCP QP.
As a distinctive feature compared to most other QP solvers, all QP formulations define a special type of variable, the slacks, which do not enter the dynamics equality constraints and which give a diagonal contribution to Hessian and inequality constraint matrices.
The exploitation of this structure allows their elimination from the QP formulation in computational complexity linear in the number of slack variables, 
making it computationally cheap to use them.
The slack variables can be used to efficiently implement soft constraints with L1 and L2 penalties.

\subsection{Dense QP}

The dense QP type describes a generic QP where Hessian and constraint matrices are assumed to be dense. 
This formulation can handle the QP sub-problems arising in single-shooting discretization schemes or in state-condensing schemes in the OCP and MPC frameworks.
\begin{align*}
\min_{v,s} & \quad \frac 1 2 \begin{bmatrix} v \\ 1 \end{bmatrix}^T \begin{bmatrix} H & g \\ g^T & 0 \end{bmatrix} \begin{bmatrix} v \\ 1 \end{bmatrix} + \\
& + \frac 1 2 \begin{bmatrix} s^{\rm l} \\ s^{\rm u} \\ 1 \end{bmatrix}^T \begin{bmatrix} Z^{\rm l} & 0 & z^{\rm l} \\ 0 & Z^{\rm u} & z^{\rm u} \\ (z^{\rm l})^T & (z^{\rm u})^T & 0 \end{bmatrix} \begin{bmatrix} s^{\rm l} \\ s^{\rm u} \\ 1 \end{bmatrix} \\
{\rm s.t.} & \quad A v = b \\
& \quad \begin{bmatrix} \underline v \\ \underline d \end{bmatrix} \leq \begin{bmatrix} J^{b,v} \\ C \end{bmatrix} v + \begin{bmatrix} J^{s,v} \\ J^{s,g} \end{bmatrix} s^{\rm l} \\
& \quad \begin{bmatrix} J^{b,v} \\ C \end{bmatrix} v - \begin{bmatrix} J^{s,v} \\ J^{s,g} \end{bmatrix} s^{\rm u} \leq \begin{bmatrix} \overline v \\ \overline d \end{bmatrix} \\
& \quad s^{\rm l}\geq \underline s^{\rm l} \\
& \quad s^{\rm u}\geq \underline s^{\rm u}
\end{align*}
The primal variables comprise generic variables $v$ and slack variables $s^{\rm l}$ ($s^{\rm u}$) associated to the lower (upper) constraints.
The Hessian matrices of the slacks $Z^{\rm l}$ and $Z^{\rm u}$ are diagonal.
The matrices $J^{\cdot , \cdot}$ are made of rows from identity matrices, and are employed to select only some components in box and soft constraints.
Furthermore, the constraint matrices are the same for the upper and the lower constraints, meaning that all constraints in the formulation are two-sided.
A mask (not represented in the above dense QP formulation) can be employed to dynamically activate or deactivate the single upper and/or lower constraints.
The considerations in this paragraph apply also to the OCP and tree OCP QP.

\subsection{OCP QP}

The OCP QP type describes a QP formulation handling many common OCP and MPC problems such as linear-quadratic regulators (LQR), constrained linear MPC problems, and QP sub-problems in sequential quadratic programming (SQP) algorithms for non-linear OCP and MPC problems.
\begin{align*}
\min_{x,u,s} & \quad \sum_{n=0}^N \frac 1 2 \begin{bmatrix} u_n \\ x_n \\ 1 \end{bmatrix}^T \begin{bmatrix} R_n & S_n & r_n \\ S_n^T & Q_n & q_n \\ r_n^T & q_n^T & 0 \end{bmatrix} \begin{bmatrix} u_n \\ x_n \\ 1 \end{bmatrix} + \\
& + \frac 1 2 \begin{bmatrix} s^{\rm l}_n \\ s^{\rm u}_n \\ 1 \end{bmatrix}^T \begin{bmatrix} Z^{\rm l}_n & 0 & z^{\rm l}_n \\ 0 & Z^{\rm u}_n & z^{\rm u}_n \\ (z^{\rm l}_n)^T & (z^{\rm u}_n)^T & 0 \end{bmatrix} \begin{bmatrix} s^{\rm l}_n \\ s^{\rm u}_n \\ 1 \end{bmatrix} \\
{\rm s.t.} & \quad x_{n+1} = A_n x_n + B_n u_n + b_n \quad , \quad n \in \mathcal H \backslash \{N\} \\ 
& \quad \begin{bmatrix} \underline u_n \\ \underline x_n \\ \underline d_n \end{bmatrix} \leq \begin{bmatrix} J_n^{b,u} & 0 \\ 0 & J_n^{b,x} \\ D_n & C_n \end{bmatrix} \begin{bmatrix} u_n \\ x_n \end{bmatrix} + \begin{bmatrix} J_n^{s,u} \\ J_n^{s,x} \\ J_n^{s,g} \end{bmatrix} s^{\rm l}_n \; , \; n \in \mathcal H \\ 
& \quad \begin{bmatrix} J_n^{b,u} & 0 \\ 0 & J_n^{b,x} \\ D_n & C_n \end{bmatrix} \begin{bmatrix} u_n \\ x_n \end{bmatrix} - \begin{bmatrix} J_n^{s,u} \\ J_n^{s,x} \\ J_n^{s,g} \end{bmatrix} s^{\rm u}_n \leq \begin{bmatrix} \overline u_n \\ \overline x_n \\ \overline d_n \end{bmatrix} \; , \; n \in \mathcal H \\ 
& \quad s_n^{\rm l}\geq \underline s_n^{\rm l} \; , \; n \in \mathcal H \\ 
& \quad s_n^{\rm u}\geq \underline s_n^{\rm u} \; , \; n \in \mathcal H 
\end{align*}
where $\mathcal H = \{0,1,\dots,N\}$.
This problem has a multi-stage structure, with cost and inequality constraints defined stage-wise, and with dynamics equality constraints coupling pairs of consecutive stages.
The primal variables are divided into state variables $x_n$, control (or input) variables $u_n$ and slack variables associated to the lower (upper) constraints $s^{\rm l}_n$ ($s^{\rm u}_n$).
The size of all variables (number of states $n_{x_n}$, number of controls $n_{u_n}$ and number of slacks $n_{s_n}$), as well as the number of box constraints $n_{b_n}$ and general polytopic constraints $n_{g_n}$ can freely vary stage-wise.
All data matrices and vectors can vary stage-wise, as required e.g. to solve QP sub-problems arising in multiple-shooting discretization schemes \citep{Bock1984}.

Note that the current formulation does not support equality constraints other than the dynamics equations, and therefore other types of equality constraints (as e.g. the constraint on the initial state value $x_0 = \hat x_0$) have to be reformulated as inequality constraints with equal upper and lower limits.
The IPM variant employed in HPIPM (see Section~\ref{sec:impl:ipm}) is able to handle such reformulation well.
This also applies to the tree OCP QP.

\subsection{Tree OCP QP}

The tree OCP QP type can handle many common robust and scenario-based OCP and MPC problems, see~\citep{Kouzoupis2019} and references therein.
\begin{align*}
\min_{x,u,s} & \quad \sum_{n \in \mathcal N} \frac 1 2 \begin{bmatrix} u_n \\ x_n \\ 1 \end{bmatrix}^T \begin{bmatrix} R_n & S_n & r_n \\ S_n^T & Q_n & q_n \\ r_n^T & q_n^T & 0 \end{bmatrix} \begin{bmatrix} u_n \\ x_n \\ 1 \end{bmatrix} + \\
& + \frac 1 2 \begin{bmatrix} s^{\rm l}_n \\ s^{\rm u}_n \\ 1 \end{bmatrix}^T \begin{bmatrix} Z^{\rm l}_n & 0 & z^{\rm l}_n \\ 0 & Z^{\rm u}_n & z^{\rm u}_n \\ (z^{\rm l}_n)^T & (z^{\rm u}_n)^T & 0 \end{bmatrix} \begin{bmatrix} s^{\rm l}_n \\ s^{\rm u}_n \\ 1 \end{bmatrix} \\
{\rm s.t.} & \quad x_m = A_m x_n + B_m u_n + b_m \quad , \quad \left\{ \begin{matrix} n \in \mathcal N \backslash \mathcal L \\ m \in \mathcal C(n) \end{matrix} \right. \\
& \quad \begin{bmatrix} \underline u_n \\ \underline x_n \\ \underline d_n \end{bmatrix} \leq \begin{bmatrix} J_n^{b,u} & 0 \\ 0 & J_n^{b,x} \\ D_n & C_n \end{bmatrix} \begin{bmatrix} u_n \\ x_n \end{bmatrix} + \begin{bmatrix} J_n^{s,u} \\ J_n^{s,x} \\ J_n^{s,g} \end{bmatrix} s^{\rm l}_n \; , \; n \in \mathcal N \\
& \quad \begin{bmatrix} J_n^{b,u} & 0 \\ 0 & J_n^{b,x} \\ D_n & C_n \end{bmatrix} \begin{bmatrix} u_n \\ x_n \end{bmatrix} - \begin{bmatrix} J_n^{s,u} \\ J_n^{s,x} \\ J_n^{s,g} \end{bmatrix} s^{\rm u}_n \leq \begin{bmatrix} \overline u_n \\ \overline x_n \\ \overline d_n \end{bmatrix} \; , \; n \in \mathcal N \\
& \quad s_n^{\rm l}\geq \underline s_n^{\rm l} \; , \; n \in \mathcal N \\
& \quad s_n^{\rm u}\geq \underline s_n^{\rm u} \; , \; n \in \mathcal N
\end{align*}
$\mathcal N$ is the set of nodes in the tree and $\hat N$ is its cardinality.
The set $\mathcal L$ contains the leaves of the tree, while $\mathcal C(n)$ denotes the set of the children of node $n$.
All data matrices and vector can vary node-wise.

\section{Algorithm}

HPIPM implements in a modular fashion all algorithmic building blocks needed to implement many variants of a primal-dual IPM.
This algorithmic richness can be used e.g. to trade off computational speed and robustness, to adapt it to the specific application.

\subsection{Primal-dual IPM}

This section contains some basic notion about primal-dual IPM algorithms.
In particular, the IPM algorithm is derived from the KKT conditions, and not from the barrier methods theory.
The interested reader can find a more detailed presentations in~\citep{Wright1997,Nocedal2006}.

Let us consider the generic QP in the form
\begin{equation}
\begin{aligned}
\min_y & \quad \frac 1 2 y^T \mathcal H y + g^T y \\
s.t. & \quad \mathcal A y = b \\
& \quad \mathcal C y \geq d
\end{aligned}
\label{eq:alg:pdipm:qp}
\end{equation}
The Lagrangian function for this QP reads
\begin{equation*}
\mathcal L (v, \pi, \lambda) = \frac 1 2 y^T \mathcal H y + g^T y - \pi^T (\mathcal A y - b) - \lambda^T (\mathcal C y - d )
\end{equation*}
where $\pi$ and $\lambda$ are the Lagrange multipliers of the equality and inequality constraints, respectively.
The first order necessary KKT optimality conditions read
\begin{subequations}
\begin{align}
& \mathcal H y + g - \mathcal A^T \pi - \mathcal C^T \lambda = 0 & \label{eq:alg:pdipm:kkt:g} \\
& - \mathcal A y + b = 0 & \label{eq:alg:pdipm:kkt:b} \\
& - \mathcal C y + d + t = 0 & \label{eq:alg:pdipm:kkt:d} \\
& \lambda_i t_i = 0 & i=1\dots,n_c \label{eq:alg:pdipm:kkt:m} \\
& (\lambda, t) \geq 0 & \label{eq:alg:pdipm:kkt:s}
\end{align}
\label{eq:alg:pdipm:kkt}
\end{subequations}
where the slack variables $t = \mathcal C y - d \geq 0$ have been introduced, and $n_c$ denotes the total number of constraints.
In case the Hessian matrix $\mathcal H$ is positive definite, the QP is convex and a point satisfying the KKT conditions is the unique global minimizer.

Equations (\ref{eq:alg:pdipm:kkt:g})-(\ref{eq:alg:pdipm:kkt:m}) are a system of nonlinear equations $f(y, \pi, \lambda, t)=0$, where the only nonlinear equation is (\ref{eq:alg:pdipm:kkt:m}).
In a nut shell, a primal-dual IPM is Newton method applied to the system of equations $f_{\tau}(y, \pi, \lambda, t)=0$ with equation (\ref{eq:alg:pdipm:kkt:m}) relaxed as
\begin{equation*}
\lambda_i t_i = \tau, \quad i=1,\dots,n_c
\end{equation*}
The homotopy parameter $\tau$ is related to the barrier parameter in barrier methods, and it is shrunk toward zero as the iterations approach the solution of (\ref{eq:alg:pdipm:kkt}).
A line search procedure is used to ensure the strict satisfaction of the inequalities (\ref{eq:alg:pdipm:kkt:s}) on the sign of Lagrange multipliers and slacks of inequality constraints.

At every iteration $k$ of the Newton method, the Newton step $(\Delta y_{\rm aff}, \Delta \pi_{\rm aff}, \Delta \lambda_{\rm aff}, \Delta t_{\rm aff})$ is found solving the linear system
\begin{multline*}
\nabla f_{\tau_k}(y_k, \pi_k, \lambda_k, t_k) \begin{bmatrix} \Delta y_{\rm aff}^T & \Delta \pi_{\rm aff}^T & \Delta \lambda_{\rm aff}^T & t_{\rm aff}^T \end{bmatrix}^T = \\
= - f_{\tau_k}(y_k, \pi_k, \lambda_k, t_k)
\end{multline*}
which takes the form
\begin{multline}
\begin{bmatrix} \mathcal H & - \mathcal A^T & - \mathcal C^T & 0 \\ - \mathcal A & 0 & 0 & 0 \\ - \mathcal C & 0 & 0 & I \\ 0 & 0 & T_k & \Lambda_k \end{bmatrix}
\begin{bmatrix} \Delta y_{\rm aff} \\ \Delta \pi_{\rm aff} \\ \Delta \lambda_{\rm aff} \\ \Delta t_{\rm aff} \end{bmatrix} =  \\
= - \begin{bmatrix} \mathcal H y_k - \mathcal A^T \pi_k - \mathcal C^T \lambda_k + g \\ - \mathcal A y_k + b\\ - \mathcal C y_k + t_k + d \\ \Lambda_k T_k e - \tau_k e \end{bmatrix} \dot = - \begin{bmatrix} r_g \\ r_b \\ r_d \\ r_m \end{bmatrix}
\label{eq:alg:pdipm:delta}
\end{multline}
The diagonal matrices $\Lambda_k$ and $T_k$ have on their diagonal the elements of the vector $\lambda_k$ and $t_k$ respectively.
The function $f_{\tau_k}(y_k, \pi_k, \lambda_k, t_k)$ at the right hand side (RHS) is denoted as the residual function.

\subsection{Delta and absolute IPM formulations}

The formulation of the Newton step in (\ref{eq:alg:pdipm:delta}) is denoted as `delta formulation' in this paper, in the sense that the solution of the linear system (\ref{eq:alg:pdipm:delta}) gives the Newton step.
By exploiting the linearity of the first three equations in (\ref{eq:alg:pdipm:kkt}), and using for the RHS of the forth equation in (\ref{eq:alg:pdipm:delta}) the expression
\begin{align*}
T_k \Delta \lambda_{\rm aff} + \Lambda_k \Delta t_{\rm aff} &= T_k (\lambda_{\rm aff} - \lambda_k) + \Lambda_k (t_{\rm aff} - t_k) \\
&= T_k \lambda_{\rm aff} + \Lambda_k t_{\rm aff} - 2 \Lambda_k T_k e
\end{align*}
the linear system (\ref{eq:alg:pdipm:delta}) can be rewritten as
\begin{equation}
\begin{bmatrix} \mathcal H & - \mathcal A^T & - \mathcal C^T & 0 \\ - \mathcal A & 0 & 0 & 0 \\ - \mathcal C & 0 & 0 & I \\ 0 & 0 & T_k & \Lambda_k \end{bmatrix}
\begin{bmatrix} y_{\rm aff} \\ \pi_{\rm aff} \\ \lambda_{\rm aff} \\ t_{\rm aff} \end{bmatrix} = 
- \begin{bmatrix} g \\ b \\ d \\ - \Lambda_k T_k e - \tau_k e \end{bmatrix}
\label{eq:alg:pdipm:abs}
\end{equation}
In this paper, this formulation is denoted as `absolute formulation', in the sense that the solution of the linear system (\ref{eq:alg:pdipm:abs}) gives the next full-step Newton iterate, and the Newton step can be computed as the difference with respect to the current iterate
\begin{multline}
\begin{bmatrix} \Delta y_{\rm aff}^T & \Delta \pi_{\rm aff}^T & \Delta \lambda_{\rm aff}^T & \Delta t_{\rm aff}^T \end{bmatrix}^T = \\
= \begin{bmatrix} y_{\rm aff}^T & \pi_{\rm aff}^T & \lambda_{\rm aff}^T & t_{\rm aff}^T \end{bmatrix}^T - \begin{bmatrix} y_k^T & \pi_k^T & \lambda_k^T & t_k^T \end{bmatrix}^T
\label{eq:alg:pdipm:step}
\end{multline}
The advantage of the absolute formulation is that (\ref{eq:alg:pdipm:step}) has a computational complexity linear in the number of variables, while the evaluation of the residual function at the RHS of (\ref{eq:alg:pdipm:delta}) involves matrix-vector multiplications, which for dense matrices have a quadratic computational complexity. 
Therefore, especially in case of systems with few variables, the absolute formulation can be computationally faster.
On the other hand, as the iterate approaches the solution, the computation of the Newton step as in (\ref{eq:alg:pdipm:step}) incurs in increasingly severe cancellation errors, as the two vectors at the RHS nearly coincide. 
This effect is exacerbated by the fact that typically in late IPM iterations the conditioning of the linear system gets increasingly worse.

\subsection{KKT system solution}

The ill-conditioning arises from the fact that in HPIPM, as in most IPM solvers, the slacks $t_k$ and the Lagrange multipliers $\lambda_k$ are eliminated as a pre-processing step before the actual factorization, as the resulting linear system is in the form of the KKT system of an equality-constrained QP, and therefore efficient factorization procedures for this class of systems can be employed.

In HPIPM, the variables $\Delta t_{\rm aff}$ and $\Delta \lambda_{\rm aff}$ are eliminated, in this order, from the KKT system in (\ref{eq:alg:pdipm:delta}), and the resulting system of linear equations
\begin{multline}
\begin{bmatrix} \mathcal H + \mathcal C^T T^{-1}_k \Lambda_k \mathcal C & - \mathcal A^T \\ - \mathcal A & 0 \end{bmatrix} \begin{bmatrix} \Delta y_{\rm aff} \\ \Delta \pi_{\rm aff} \end{bmatrix} =  \\
= - \begin{bmatrix} r_g - \mathcal C^T(T^{-1}_k \Lambda_k r_d - T^{-1}_k r_m ) \\ r_b \end{bmatrix}
\label{eq:alg:kkt}
\end{multline}
which is in the form of the KKT system of an equality constrained QP, is solved with structure-exploiting routines operating on dense sub-matrices in case of all QP types currently supported.
This approach corresponds to statically selecting the pivoting sequence based solely on sparsity, and hand-crafting efficient factorization algorithms tailored to the QP type at hand.
The alternative approach, which is currently not employed in HPIPM, is to directly factorize the KKT system in (\ref{eq:alg:pdipm:delta}) using sparse LA, which 
may additionally perform dynamic pivoting for stability in an attempt to keep the conditioning under control.

In HPIPM, in the case of the dense QP type, the linear system in (\ref{eq:alg:kkt}) is solved using either the Schur-complement or the null-space methods~\citep{Nocedal2006}.
In case the dense QPs arise from the condensing of the OCP QP type, there are no equality constraints, and the factorization of the KKT matrix in (\ref{eq:alg:kkt}) can be simply performed using the dense Cholesky factorization.

In the case of the OCP QP type, as proposed in the MPC context by~\citep{Steinbach1994} and similarly to HPMPC~\citep{Frison2014}, the linear system in (\ref{eq:alg:kkt}) is factorized using a structure-exploiting backward Riccati recursion.
Two versions of the Riccati recursion are available~\citep{Frison2015a}: (a) the classical implementation (which requires the reduced Hessian with respect to the dynamics equality constraints to be positive definite, but allows the full-space Hessian to be indefinite), and (b) the square-root implementation (which in order to reduce the flop count employs the Cholesky factorization of the Riccati recursion matrix, and therefore requires the full-space Hessian to be positive definite).
In case the number of states $x_n$ and inputs $n_u$ is constant across stages, the computational complexity is of $\mathcal O(N(n_x+n_u)^3)$ flops, which is linear in the number of stages.

In the case of the tree OCP QP type, the linear system in (\ref{eq:alg:kkt}) is factorized using a Riccati recursion modified to exploit the tree structure~\citep{Frison2017a}.
This factorization proceeds from the leaves of the tree, and therefore it does not introduce any fill-in outside the problem data blocks.
The Riccati recursion can be implemented using the same two variants described for the OCP QP type.
In case the number of states $x_n$ and inputs $n_u$ is constant across stages and realizations, the computational complexity is of $\mathcal O(\hat N (n_x+n_u)^3)$, and therefore it is linear in the number of nodes in the tree $\hat N$.

\subsection{Iterative refinement}

In practice, due to the finite accuracy of numerical computations, any procedure for the solution of the KKT system (\ref{eq:alg:pdipm:delta}) can compute the solution only up to a certain accuracy.
Especially in case of unstable or ill-conditioned systems, and in late IPM iterations, the ill-conditioning of the KKT system (\ref{eq:alg:pdipm:delta}) can make the accuracy in the Newton step too low to be useful for the IPM algorithm.
In such cases, it may be beneficial to employ a few iterative refinement steps in order to attempt to improve the accuracy of the Newton step.

The iterative refinement residuals are defined as the residuals in the solution of the KKT system (\ref{eq:alg:pdipm:delta}),
\begin{equation*}
\begin{bmatrix} r_g^{\rm ir} \\ r_b^{\rm ir} \\ r_d^{\rm ir} \\ r_m^{\rm ir} \end{bmatrix} \dot =
\begin{bmatrix} \mathcal H & - \mathcal A^T & - \mathcal C^T & 0 \\ - \mathcal A & 0 & 0 & 0 \\ - \mathcal C & 0 & 0 & I \\ 0 & 0 & T_k & \Lambda_k \end{bmatrix}
\begin{bmatrix} \Delta y_{\rm aff} \\ \Delta \pi_{\rm aff} \\ \Delta \lambda_{\rm aff} \\ \Delta t_{\rm aff} \end{bmatrix} +
\begin{bmatrix} r_g \\ r_b \\ r_d \\ r_m \end{bmatrix}
\label{eq:alg:itref}
\end{equation*}
In case they are not small enough, it is possible to solve again a system in the form (\ref{eq:alg:pdipm:delta}) with the same KKT matrix at the left hand side (LHS) and the iterative refinement residuals as RHS, and compute a correction term to iteratively improve the accuracy of the Newton step.
Since no new matrix factorization is performed, the computation of this correction term is cheap compared to the factorization and solution of the KKT system in (\ref{eq:alg:pdipm:delta}), especially in case of large systems.

\subsection{Condensing and partial condensing}

HPIPM provides a condensing module, which implements routines to convert an OCP QP into a dense QP by eliminating the state variables at all stages except stage 0.
This is possible because such state variables are not real degrees of freedom, as they are determined uniquely from the inputs and the state at stage 0 by using the dynamics equality constraints.
The partial (or block) condensing module, which converts an OCP QP into another OCP QP with smaller number of stages, is built on top of the condensing module in a modular fashion, by performing condensing on blocks of stages.

There exist at least three Hessian condensing algorithms with different computational complexities, ranging between quadratic and cubic in both the number of states and the number of stages~\citep{Frison2015a}.
The algorithm currently implemented in HPIPM is the algorithm proposed in~\citep{Frison2016}, which has a computational complexity quadratic in the number of stages and cubic in the number of states.
This algorithm is based on a backward recursion which has similarities to the Riccati recursion.
As such, two versions of this condensing algorithm are available: one analogous to the classical Riccati recursion and one analogous to the square-root Riccati recursion, with similar requirements about Hessian positive definiteness.
This condensing algorithm is the fastest when the state at stage 0 is kept as an optimization variable (such as in partial condensing) 
and in case of long horizons.

\section{Implementation}


\subsection{Linear algebra}

Algorithms in HPIPM are explicitly designed to exploit the overall structure of the optimization algorithms.
Sparsity within data and working matrices is not considered.
Therefore, computations are cast in term of operations on matrices with fixed structure, like e.g. dense, symmetric, triangular or diagonal.

In HPIPM, the KKT modules contain the most computationally expensive routines, such as KKT system factorization and solution, and residuals computation.
Such operations have computational complexities cubic or quadratic in the number of (stage) variables.
In case of all currently implemented QP types, in the KKT modules HPIPM employs the structure-based interface in the high-performance LA library BLASFEO~\citep{Frison2018} for such operations.
These routines are designed to give close-to-peak performance for matrices of moderate size (assumed to fit in cache) and are optimized for several hardware architectures (e.g. exploiting different classes of SIMD instructions).

The IPM core operations in HPIPM are independent of the specific QP type and solely operate on unstructured arrays of doubles holding e.g. the set of all primal variables, or Lagrange multipliers or slack variables.
These operations have computational complexity linear in the number of variables, and are currently offered in two implementations, one coded in generic C and one using AVX intrinsics. 

\subsection{IPM implementation choices} \label{sec:impl:ipm}

At a high-level, HPIPM implements an infeasible-start Mehrotra's predictor-corrector IPM~\citep{Mehrotra1992}, proposed in the MPC context e.g. in~\citep{Rao1998}.
This standard algorithmic scheme is tweaked by several arguments, whose different choice gives rise to a family of algorithms with different trade-offs between computational speed and reliability.
A `mode' argument is introduced in order to select pre-defined sets of arguments tailored to a particular purpose.
The modes currently available are:
\begin{description}
\item[speed\_abs]
This mode selects an algorithm analogous to the one implemented in HPMPC.
The focus lies on the highest solution speed, and the IPM is based on the absolute formulation.
The residuals are optionally computed only before returning, and therefore the IPM loop exit conditions are only based on duality measure tolerance, maximum number of iterations and minimum step length.
This algorithm is indicated for small and well-conditioned systems, in case tight bounds on the solution accuracy are not required.
\item[speed]
The focus lies on the high solution speed, but the more reliable IPM delta formulation is employed.
The KKT residuals are computed at each IPM iteration, and exit conditions are additionally based on minimum residuals tolerance.
\item[balance]
In this mode, the accuracy of the KKT factorization is checked by computing its iterative refinement residuals.
In case accuracy is too low, the factorization is repeated by replacing all Cholesky factorizations of normal matrices in the form $A\cdot A^T$ with array algorithms based on QR factorizations, which have better numerical properties as they never explicitly form the worse-conditioned normal matrix $A\cdot A^T$.
Furthermore, a few iterative refinement steps are performed on the prediction-centering-correction direction, in case this is not accurate enough.
\item[robust]
In this mode, the KKT factorization is always performed using the more accurate array algorithms based on QR factorizations.
A higher number of iterative refinement steps is performed on the prediction-centering-correction direction, in case this is not accurate enough.
\end{description}

Some additional implementation choices are common to all modes.
The Mehrotra's corrector step is implemented conditionally, and is only applied if it is not increasing the duality measure too much compared to the affine step; otherwise, a predictor-centering step is applied.
The Hessian matrix of both primal and dual variables can be regularized by adding a scaled identity matrix during factorization; this regularization does not affect the original problem data, and therefore the effect of regularization on the solution can be offset by means of iterative refinement.
The value of the inequality constraint Lagrange multipliers and slacks is bounded to a minimum value, in order to give the possibility to bound the ill-conditioning of the KKT system in late IPM iterations.
Finally, it is possible to provide a guess to warm start either the primal, or both the primal and the dual variables (this second option is particularly effective in SQP algorithms, since close to the NLP solution a good guess for the primal variables is zero, and for the dual variables coincides with the Lagrange multipliers of the NLP).

\section{Numerical experiments}

In this section, numerical experiments evaluate the relative speed (Section~\ref{sec:exp:lmss}) and robustness (Section~\ref{sec:exp:cuter}) of the different HPIPM modes, and compare HPIPM to other state-of-the-art MPC solvers (Section~\ref{sec:exp:chain}).


\subsection{Linear mass spring system} \label{sec:exp:lmss}

The linear mass spring system is a scalable benchmark commonly employed to evaluate linear MPC solvers~\citep{Wang2010a}.
The variant proposed in~\citep{Domahidi2012} is considered, comprising a horizontal chain of $M$ masses connected to each other and to walls using springs, and where the first $M-1$ masses are controlled using a force.
Therefore the number of states $n_x=2M$, inputs $n_u=M-1$ and box constraints $n_b=n_x+n_u=3M-1$ per stage scales with the number of masses.

Table~\ref{tab:lmpc} compares the solution time for 3 HPIPM modes with HPMPC.
The algorithm in HPIPM with {\tt  speed\_abs} mode is analogous to the one in HPMPC, and the computational times are also within 10-15\% of each other; HPIPM is slightly slower due to its more modular nature, and its use of an external LA library.
The HPIPM {\tt speed} mode is about 2 times slower for small system, but the difference decreases for large systems (since the same KKT factorization is employed).
The HPIPM {\tt robust} mode is roughly 3 times slower for all systems, since a slower but more accurate KKT factorization is employed.

\begin{table}
\caption{Comparison of solvers for box-constrained linear MPC: {\tt speed\_abs}, {\tt speed} and {\tt robust} modes for HPIPM, low-level interfaces for HPMPC. 
$M$ is the number of masses; 
$N$ is the horizon length.
Runtimes in seconds; number of IPM iterations fixed to 10.
Intel Core i7 4800MQ CPU with `Haswell' target in BLASFEO and `AVX2' target in HPMPC.}
\centering
{\begin{tabular}{cc|lll|l}
\toprule
& & HPIPM & HPIPM & HPIPM & HPMPC \\ 
$M$ & $N$ & speed\_abs & speed & robust & low-level \\ 
\midrule
 2 & 10 &  $6.71\cdot10^{-5}$ &  $1.17\cdot10^{-4}$ &  $1.80\cdot10^{-4}$ &  $5.7\cdot10^{-5}$  \\ 
 4 & 10 &  $1.10\cdot10^{-4}$ &  $1.82\cdot10^{-4}$ &  $3.67\cdot10^{-4}$ &  $8.9\cdot10^{-5}$  \\ 
 6 & 30 &  $6.36\cdot10^{-4}$ &  $9.00\cdot10^{-4}$ &  $1.96\cdot10^{-3}$ &  $4.88\cdot10^{-4}$ \\ 
11 & 10 &  $4.62\cdot10^{-4}$ &  $6.22\cdot10^{-4}$ &  $1.61\cdot10^{-3}$ &  $4.05\cdot10^{-4}$ \\ 
15 & 10 &  $8.29\cdot10^{-4}$ &  $1.05\cdot10^{-3}$ &  $2.94\cdot10^{-3}$ &  $7.51\cdot10^{-4}$ \\ 
30 & 30 &  $1.21\cdot10^{-2}$ &  $1.45\cdot10^{-2}$ &  $4.53\cdot10^{-2}$ &  $1.17\cdot10^{-2}$ \\ 
\bottomrule
\end{tabular}}
\label{tab:lmpc}
\end{table}

\subsection{CUTEr} \label{sec:exp:cuter}

This section evaluates the performance and robustness of HPIPM using the challenging Maros-M\'esz\'aros test set of convex QPs~\citep{Maros1999}, which can be accessed though the CUTEr testing environment.
In particular, the QPs available on the qpOASES COIN-OR repository are employed, which are the 43 problems with up to 250 variables and up to 1001 constraints, as discussed in~\citep{Ferreau2014}.
All matrices are in dense format.

\begin{figure}[tpb]
\centering
\includegraphics[width=0.9\linewidth]{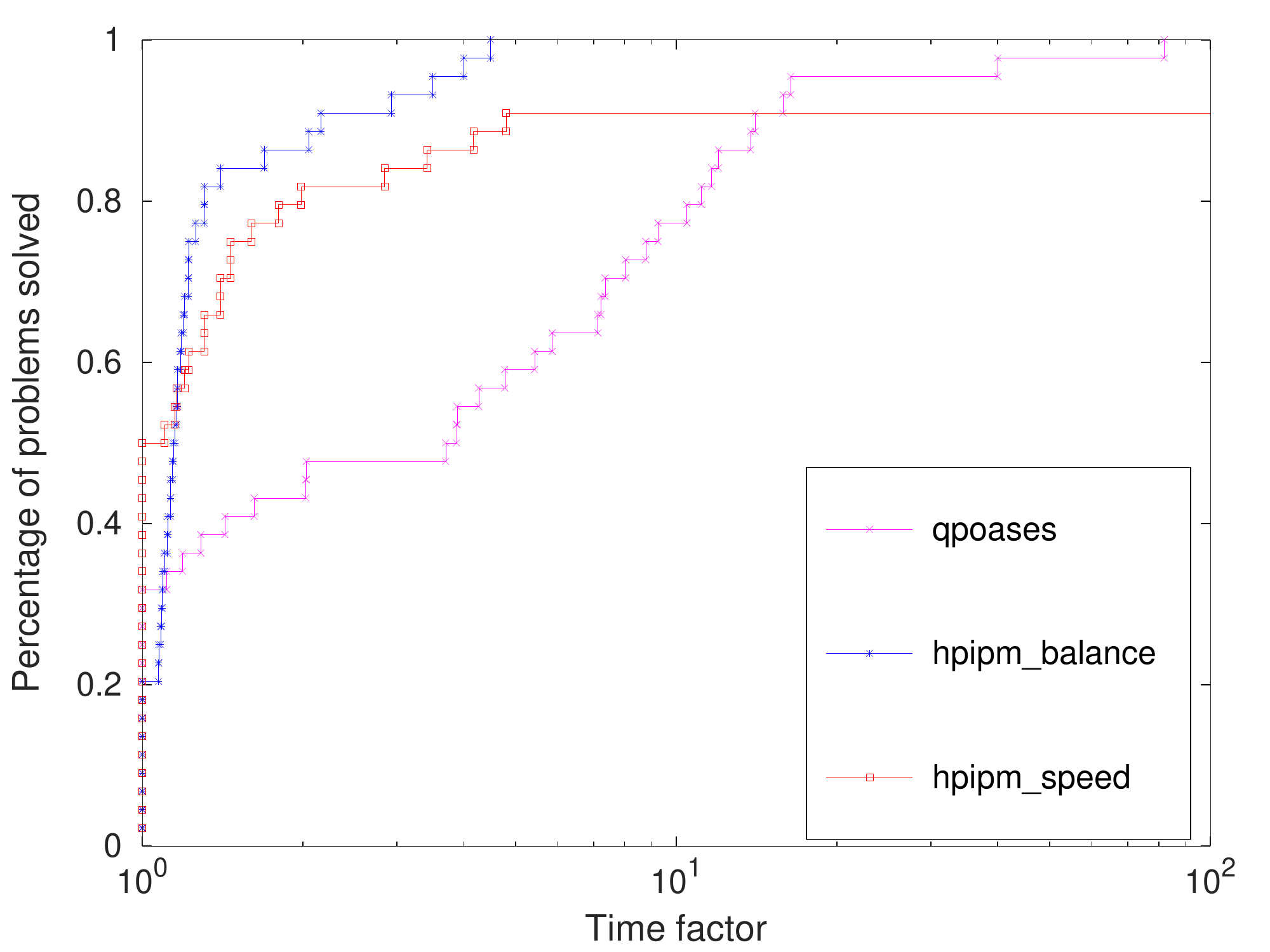}
\caption{Performance graph on the 43 problems from the Maros-M\'esz\'aros QP test set with up to 250 variables and up to 1001 constraints.
Intel Core i7 4810MQ CPU with `Haswell' target in BLASFEO.}
\label{fig:perf}
\end{figure}

Figure~\ref{fig:perf} contains a performance graph comparing the performance of the {\tt speed} and {\tt balance} modes in HPIPM with qpOASES; the required accuracy is set to $10^{-6}$ for the infinity norm of residuals.
HPIPM in {\tt speed} mode is the fastest solver in about 50\% of the problems, but it fails to solve 4 problems. 
HPIPM in {\tt balance} mode can solve all problems, it is the fastest solver in about 20\% -- and within a factor 2 of the fastest solver for more than 90\% -- of the problems. 
qpOASES is the fastest solver in the remaining 30\% of the problems, but it gets quickly outperformed (often by a factor 10 or more) in the remaining problems.

HPIPM in {\tt speed\_abs} mode (not shown in Figure~\ref{fig:perf}) would fail to converge in 32 problems, and return an inaccurate solution in the remaining problems.
This shows the great improvement in reliability over HPMPC (which implements an algorithm analogous to HPIPM in {\tt speed\_abs} mode, but can not be directly compared as it does not provide a dense QP solver) in case of challenging QPs.

\subsection{Chain of masses} \label{sec:exp:chain}

This section compares the performance of HPIPM against other state-of-the-art QP solvers for MPC problems.
In particular, the NLP framework acados is employed to generate constrained time-varying linear MPC problems.
The test problem is the nonlinear hanging chain problem used for the comparison of QP solvers in~\citep{Kouzoupis2018}, where the chain of masses starts at rest position, and a strong perturbation is generated by overwriting the controller output in the first 5 simulation steps; afterwards the controller takes over again and takes the chain of masses back to the resting position.
Four free masses are employed, accounting for a total of 24 states (the 3D position and velocity of each free mass), 3 controls (the 3D force acting on the last mass), 4 state and 3 control bounds; the control horizon is of 8 seconds, divided into 40 equidistant stages. 

\begin{figure}[tpb]
\centering
\includegraphics[width=0.88\linewidth]{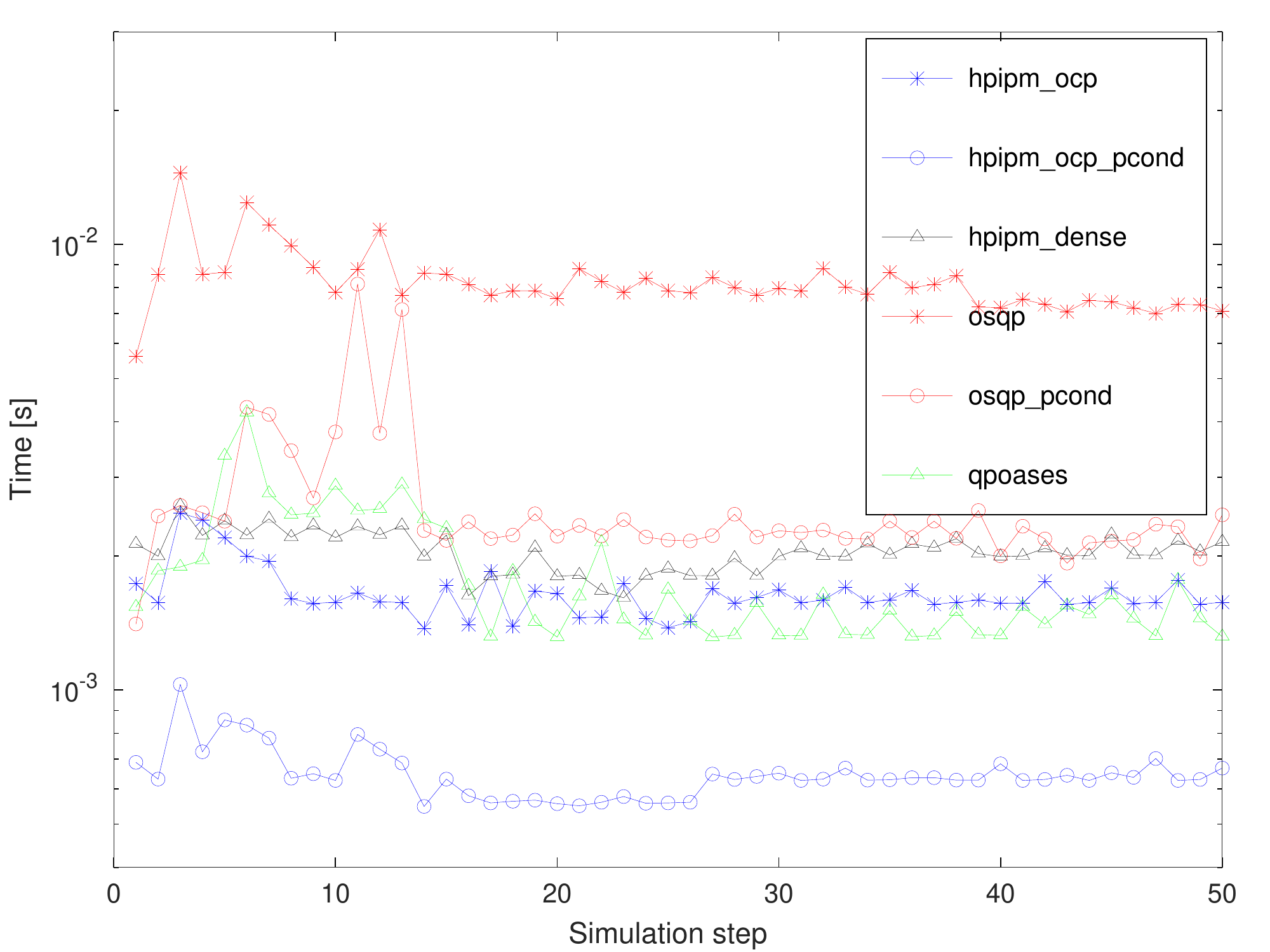}
\caption{Comparison of 
QP solvers in the solution of QPs generated using acados in a closed loop MPC simulation on the nonlinear chain of masses.
Intel Core i7 4810MQ CPU with `Haswell' target in BLASFEO.}
\label{fig:chain}
\end{figure}

The numerical experiment consists of a closed loop simulation of 50 steps, where the real time iteration scheme in acados is used to generate one QP to be solved per simulation step.
Optionally partial condensing is employed to reformulate the OCP QP into another OCP QP with 5 stages.
The QP solution times are in Figure~\ref{fig:chain} for the solvers HPIPM ({\tt speed} mode; OCP QP solver for the fully sparse or partially condensed QP; dense QP solver for the condensed QP), OSQP (sparse solver for the fully sparse or partially condensed QP) and qpOASES (for the condensed QP).
(Partial) condensing is provided by the acados wrapper to the corresponding HPIPM routines.
Note that only the call to the QP solver is timed: therefore, the (partial) condensing time is not included, nor the SQP algorithm time.
Warm start is employed for OSQP and qpOASES, as these solvers are found to benefit from it; OSQP employs a polishing step to improve the low accuracy of the ADMM (Alternating Direction Method of Multipliers) solution.

HPIPM solving the partially condensed OCP QP stands out as the fastest solver, with solution times steadily below 1 millisecond.
Both HPIPM and OSQP are significantly faster in solving the partially condensed QP compared to the corresponding fully sparse QP.
In general, HPIPM (being an IPM) shows a rather constant solution time, while qpOASES (active set) and OSQP (ADMM) require more time (due to more iterations) when the QP problems are more challenging, at the beginning of the simulation.

\section{Conclusion}

This paper introduced HPIPM, a high-performance QP framework for the solution of MPC problems.
HPIPM provides IPM solvers for dense, OCP and tree OCP QPs, as well as (partial) condensing routines.
In particular, compared to the predecessor HPMPC solver, reliability is shown to be significantly improved in case of challenging QPs.
Numerical experiments show that HPIPM speed excels against other state-of-the-art QP solvers for MPC. 


\end{document}